\documentclass[11pt]{amsart}

\usepackage{amssymb}
\usepackage{latexsym}
\usepackage{amsmath}
\usepackage{amsthm}
\usepackage{amsfonts}
\usepackage{color}
\usepackage{pdfsync}
\usepackage[usenames,dvipsnames]{pstricks}
\usepackage{epsfig}
\usepackage{pst-grad} 
\usepackage{pst-plot} 

\numberwithin{equation}{section}

\usepackage{amssymb}

\newcommand{\beq}{\begin{equation} }
\newcommand{\eqq}{\end{equation} }
\newcommand{\cuad}{{\sqcap\kern-.68em\sqcup}}

\newtheorem{remark}{Remark}[section]
\newcommand{\bremark}{\begin{remark} \em}
\newcommand{\eremark}{\end{remark} }

\def\beeq{\begin{equation}}
\def\eeq{\end{equation}}
\newcommand{\begeqaet}{\begin{eqnarray*}}
\newcommand{\eneqaet}{\end{eqnarray*}}

\let\Section=\section
\def\section{\setcounter{equation}{0}\Section}
\newtheorem{Lem}{Lemma}[section]
\newtheorem{Thm}{Theorem}[section]

\begin{document}
\begin{center}{\bf\Large Existence and symmetry result for Fractional p-Laplacian in $\mathbb{R}^{n}$  }\medskip

\bigskip

\bigskip

{C\'esar E. Torres Ledesma}

 Departamento de Matem\'atica ,\\
 Universidad Nacional de Trujillo\\
Av. Juan Pablo II s/n, Trujillo, Per\'u.\\
 {\sl  (ctl\_576@yahoo.es)}


\medskip

\medskip
\medskip
\medskip
\medskip

\end{center}

\centerline{\bf Abstract}

\medskip

In this article we are interested in  the following fractional $p$-Laplacian equation in $\mathbb{R}^n$
\begin{eqnarray*}
&(-\Delta)_{p}^{\alpha}u + V(x)u^{p-2}u = f(x,u)   \mbox{ in } \mathbb{R}^{n},   
\end{eqnarray*}
where $p\geq 2$, $0< s < 1$, $n\geq 2$ and subcritical p-superlinear nonlinearity. By using mountain pass theorem with Cerami condition we prove the existence of nontrivial solution. Furthermore, we show that this solution is radially simmetry. 

\medskip

\date{}

\setcounter{equation}{0}
\section{ Introduction}

Recently, a great attention has been focused on the study of problem involving fractional and non-local operators. This type of problem arises in many different applications, such as, continuum mechanics, phase transition phenomena, population dynamics and game theory, as they are the typical outcome of stochastically stabilization of L\'evy processes, see \cite{DA, LC, NL-1, RMJK} and the references therein. The literature on non-local operators and their applications is very interesting and quite large, we refer the interested reader to \cite{LAGPLM, BBECAP, LCLS, LCJROS, AIMS-1, AIMS-2, AISLKPMS, PFAQJT, PFCT-1, PFCT-2, KT} and the references therein. For the basic properties of fractional Sobolev spaces, we refer the interested reader to \cite{EDNGPEV}. 
 
 In this paper, our aim is to show the existence of weak solution for the following class of equations:
\begin{equation}\label{00}
(-\Delta)_{p}^{s}u + V(x)u^{p-2}u = f(x,u)   \mbox{ in } \mathbb{R}^{n},
\end{equation}
where $p\geq 2$, $0< s < 1$, $n\geq 2$ and $(-\Delta)_{p}^{s}$ is the fractional p-laplacian defined by
\begin{equation}\label{01}
(-\Delta)_{p}^{2}u(x) = 2\lim_{\epsilon \searrow 0} \int_{\mathbb{R}^{n} \setminus B_{\epsilon}(x)} \frac{|u(x) - u(y)|^{p-2}(u(x) - u(y))}{|x-y|^{p-2}}dy,\;\;x\in \mathbb{R}^{n}.
\end{equation}

\noindent
In order to obtain the existence of weak solution for (\ref{00}), let us recall some result related to the fractional Sobolev space $W^{s,p}(\mathbb{R}^n)$. First of all, define the Gagliardo seminorm by
$$
[u]_{s,p} = \left(\int_{\mathbb{R}^{n}}\int_{\mathbb{R}^{n}} \frac{|u(x) - u(y)|^p}{|x-y|^{n+sp}}dxdy \right)^{1/p},
$$
where $u: \mathbb{R}^{n} \to \mathbb{R}$ is a measurable function. Now, consider that the fractional Sobolev space given by
$$
W^{s,p}(\mathbb{R}^n) = \{u\in L^{p}(\mathbb{R}^{n})\;:\;u\; \mbox{is measurable and }\;[u]_{s,p}<\infty\}
$$
is assumed to be endowed with norm
$$
\|u\|_{s,p} = \left([u]_{s,p} + \|u\|_{p}^{p} \right)^{1/p},
$$
where the fractional critical exponent is defined by
$$
p_{s}^{*} = \left\{ \begin{array}{ll}
\frac{np}{n-sp},&\mbox{if},\;sp<n;\\
\infty,&\mbox{if}\;sp\geq n
\end{array}
\right.
$$

Moreover we consider the fractional Sobolev space with potential
$$
X^{s}:=\left\{ u\in W^{s,p}(\mathbb{R}^n):\;\;\int_{\mathbb{R}^n} V(x)|u|^pdx < \infty \right\}
$$
endowed with the norm
$$
\|u\|_{X^s} = \left([u]_{s,p}^{p} + \|V(x)^{1/p}u\|_{p}^{p} \right)^{1/p}.
$$
Now, we make the following assumptions on the functions $V$ and $f$. 
\begin{enumerate}
\item[($V$)] $V\in C(\mathbb{R}^{n})$, $\inf_{\mathbb{R}^{n}}V(x) \geq V_{0}>0$ and $\mu (x\in \mathbb{R}^n:V(x) \leq M) < +\infty$, $\forall M>0$.
\item[$(f_{1})$] $f\in C(\mathbb{R}^{n} \times \mathbb{R})$ and satisfies 
\begin{eqnarray}\label{03}
\lim_{|t|\to \infty} \frac{f(x,t)}{|t|^{q-1} } = 0,\;\lim_{|t|\to \infty} \frac{f(x,t)t}{|t|^{p}} = +\infty
\end{eqnarray}
uniformly in $x\in \mathbb{R}^n$ for some $q\in (p,p_{s}^{*})$, where $F(x,t) = \int_{0}^{t}f(x,s)ds$.
\item[$(f_{2})$] $f(x,t) = o(|t|^{p-2}t)$ as $|t|\to 0$, uniformly in $x\in \mathbb{R}^n$.
\item[$(f_{3})$] There exists $\theta \geq 1$ such that $\theta \mathcal{F}(x,t) \geq \mathcal{F}(x,\sigma t)$ for $(x,t) \in \mathbb{R}^n \times \mathbb{R}$ and $\sigma \in [0,1]$, where
$$
\mathcal{F}(x,t) = f(x,t)t - pF(x,t).
$$
\end{enumerate}
By condition $(f_{1})$ and ($f_{2}$), for any $\epsilon >0$ there exists a constant $C_\epsilon>0$ such that 
\begin{equation}\label{04}
|F(x,t)| \leq \frac{\epsilon}{p} |t|^{p} + \frac{C_{\epsilon}}{q}|t|^{q}.
\end{equation}
Consequently, the energy functional $I: X^{s}(\mathbb{R}^n) \to \mathbb{R}$,
\begin{equation}\label{03}
I(u) = \frac{1}{p} \|u\|_{X^s}^{p} - \int_{\mathbb{R}^n}F(x,u)dx
\end{equation}
is well defined and of class $C^1$. The derivative of $I$ is given by 
\begin{eqnarray*}
I'(u)v &=&\int_{\mathbb{R}^n}\!\!\int_{\mathbb{R}^n} \frac{|u(x) - u(y)|^{p-2}(u(x)-u(y))(v(x) - v(y))}{|x-y|^{n+sp}}dxdy +\! \int_{\mathbb{R}^n}V(x)uvdx\nonumber\\
&&- \int_{\mathbb{R}^n}f(x,u)vdx.
\end{eqnarray*}
for $v\in X^s$. Therefore, the critical points of $I$ are weak solutions of (\ref{00}). Now we are ready to state our main result

\begin{Thm}\label{Itm}
Suppose that the conditions $(V),(f_1)-(f_3)$ hold. Then, the problem (\ref{00}) has a nontrivial solution.
\end{Thm}

In our second main Theorem we consider the question of the symmetry properties for the solution obtained by Theorem \ref{Itm}. First, we mention some previous result obtained in the case $p=2$. Dipierro, Patalucci and Valdinoci \cite{SDGPEV}, consider the existence of radially symmetric solutions of (\ref{00}) when $V$ and $f$ do not depend explicitly on the space variable $x$. For the first time, using rearrangement tools and following the ideas of Berestycki and Lions \cite{HBPL}, the authors prove existence of a nontrivial, radially symmetric, solution to
\begin{eqnarray}\label{Ieq05}
(-\Delta)^{\alpha}u + u = |u|^{p-1}u &\mbox{in } \mathbb{R}^{n},\;\;u \in H^{\alpha}(\mathbb{R}^{n}).
\end{eqnarray}
Solutions of (\ref{Ieq05}) can be obtained by finding critical points of the Euler-Lagrange functional $I$ defined in the fractional Sobolev spaces $H^{\alpha}(\mathbb{R}^{n})$ by
$$
I(u) = \frac{1}{2}\int_{\mathbb{R}^{n}}\int_{\mathbb{R}^{n}} \frac{|u(x) - u(z)|^{2}}{|x-z|^{n+2\alpha}}dzdx + \int_{\mathbb{R}^{n}} \left( \frac{1}{2}|u(x)|^{2} - \frac{1}{p+1}|u(x)|^{p+1} \right)dx.
$$
They noted that, by the fractional Polya-Szeg\"o inequality
\begin{equation}\label{Ieq06}
I(u^{*}) \leq I(u)
\end{equation}
where $u^{*}$ is the symmetric rearrangement of $u$. Therefore if there is a minimizer of $I$, it must be a symmetric minimizer. Their proof is based on a minimization method, working with an appropriate constraint. This constraint is useful in this case because of the autonomous character of (\ref{Ieq05}) and the fact that one can rescale in $\mathbb{R}^{n}$. Felmer, Quaas and Tan \cite{PFAQJT}, studied symmetry of positive solution of the nonlinear fractional Schr\"odinger equation
$$
(-\Delta)^{\alpha}u + u = f(u) \;\;\mbox{in}\;\;\mathbb{R}^{n},
$$
using the integral form of the moving planes method, assuming that the function $f$ is super-linear, with sub-critical growth at infinity and
 \begin{itemize}
 \item[(F)] $f\in C^{1}(\mathbb{R})$, increasing and there exists $\tau >0$ such that
 $$
 \lim_{v \to 0} \frac{f'(v)}{v^{\tau}} = 0,
 $$
 \end{itemize}
 taking advantage of the representation formula for $u$ given by 
\begin{equation}\label{Ieq07}
u(x) = (\mathcal{K} * f(u))(x),\;\;x \in \mathbb{R}^{n},
\end{equation}
where the kernel $\mathcal{K}$, associated to the linear part of the equation, plays a key role in the moving planes argument. 

Very recently Felmer and Torres \cite{PFCT-1, PFCT-2}, considered positive solutions of nonlinear Schr\"odinger equation with non-local regional diffusion 
\begin{equation}\label{Eq04}
\epsilon^{2\alpha}(-\Delta)_{\rho}^{\alpha}u + u = f(u)\;\;\mbox{in}\;\;\mathbb{R}^n,\;\;u\in H^{\alpha}(\mathbb{R}^n).
\end{equation}  
The operator $(-\Delta)_{\rho}^{\alpha}$ is a variational version of the non-local regional Laplacian, defined by
$$
\int_{\mathbb{R}^n}(-\Delta)_{\rho}^{\alpha}uvdx = \int_{\mathbb{R}^n}\int_{B(0,\rho (x))} \frac{[u(x+z) - u(x)][v(x+z)-v(x)]}{|z|^{n+2\alpha}}dzdx.
$$ 
Under suitable assumptions on the nonlinearity $f$ and the range of scope $\rho$, they obtained the existence of a ground state by mountain pass argument and a comparison method. Furthermore, they analyzed symmetry properties and concentration phenomena of these solutions. These regional operators present various interesting characteristics that make them very attractive from the point of view of mathematical theory of non-local operators.  

Inspired by these previous results, we consider that the nonlinearity $f(x,u) = f(u)$ satisfies ($f_1$)-($f_3$) and  regarding the potential $V$ we assume 
\begin{itemize}
\item[($V_{r}$)] $V$ is radially symmetric and increasing.
\end{itemize}

Now we state the main Theorem in our paper.

\begin{Thm}\label{SRtm1}
Suppose that $(f_{1})-(f_{3})$ and $(V)-(V_{r})$ hold. Then the mountain pass value is achieved by a radially symmetric function, which  is a  weak solution of (\ref{00}).  
\end{Thm}


To prove this theorem we follow the ideas of \cite{PFCT-2}, namely, we proceed by using rearrangements and variational methods. The idea to prove our result, consists in replacing the path $\gamma$ in the mountain pass setting by its symmetrization $\gamma^{*}: t \in [0,1] \to \gamma (t)^{*}$. Then $u$ would be near of the set $\gamma^{*}([0,1])$. This idea works since rearrangements are continuous in $H^{\alpha}(\mathbb{R}^{n})$ see \cite{FAEL}.

We note that the approach used in \cite{PFAQJT} is not possible to be used for problem (\ref{00}), since a representation formula like (\ref{Ieq07}) is not available in general for $(-\Delta)_{p}^{s}$. On the other hand we cannot use the approach in \cite{SDGPEV} since our problem is $x$-dependent and we can not use scale changes in $\mathbb{R}^{n}$.


\section{Preliminaries and notation}

To study the fractional problem (\ref{00}) the so-called fractional Sobolev spaces $W^{s,p}(\mathbb{R}^n)$ with $0<s<1$ are expedient. If $1<p<\infty$, as usual, the norm is defined through 
$$
\|u\|_{s,p}^{p} = \int_{\mathbb{R}^n}\int_{\mathbb{R}^n} \frac{|u(x) - u(y)|^p}{|x-y|^{n+sp}}dxdy + \int_{\mathbb{R}^n}|u(x)|^pdx.
$$
We recall the Sobolev embedding theorem.
\begin{Thm}\label{X1tm}
\cite{EDNGPEV} Let $s\in (0,1)$ and $p\in [1,+\infty)$ be such that $sp<n$. Then there exists a positive constant $C = C(n,p,s)$ such that
$$
\|u\|_{L^{p_{s}^{*}}}^{p} \leq C \int_{\mathbb{R}^n}\int_{\mathbb{R}^n} \frac{|u(x) - u(y)|^p}{|x-y|^{n+sp}}dxdy, 
$$
where $p_{s}^{*} = \frac{np}{n-sp}$ is the so-called ``fractional critical exponent''. Consequently, the space $W^{s,p}(\mathbb{R}^n)$ is continuously embedded in $L^{q}(\mathbb{R}^n)$ for any $q\in [p,p_{s}^{*}]$. Moreover the embedding $W^{s,p}(\mathbb{R}^n)  \hookrightarrow  L_{loc}^{q}(\mathbb{R}^n)$ is compact for $q\in [p, p_{s}^{*})$.
\end{Thm}

Now consider the space $X^{s}$ defined by
$$
X^{s}:=\left\{ u\in W^{s,p}(\mathbb{R}^n):\;\;\int_{\mathbb{R}^n} V(x)|u|^pdx < \infty \right\},
$$
endowed with the norm
$$
\|u\|_{X^s} = \left(\iint_{\mathbb{R}^n \times \mathbb{R}^n}\frac{|u(x) + u(y)|^p}{|x-y|^{n+sp}}dxdy + \int_{\mathbb{R}^n}V(x)|u(x)|^pdx \right)^{1/p}.
$$
From $(V_{1})$, Theorem \ref{X1tm} and H\"older inequality, we have $X^s \hookrightarrow L^{q}(\mathbb{R}^n)$ for $p\leq q \leq p_{s}^{*}$. Moreover, the following compactness result holds. It was proved in \cite{XC} in the case $p=2$. For the general case, the proof is similar. We give it here for reader's convenience.  

\begin{Lem}\label{X1lem}
Suppose that $(V_{1})$ hold. Then $X^s  \hookrightarrow L^{q}(\mathbb{R}^n)$ is compact for $q\in [p, p_{s}^{*})$
\end{Lem}

\noindent 
{\bf Proof.} Let $\{u_{n}\} \in X^s$ be a bounded sequence of $X^s$ such that $u_{n} \rightharpoonup 0$ in $X^s$. Then, by Theorem \ref{X1tm}, $u_{n} \to 0$ in $L_{loc}^{q}(\mathbb{R}^n)$ for $p\leq q < p_{s}^{*}$. We claim that
\begin{equation}\label{s01}
u_{n} \to 0\;\;\mbox{strongly in}\;\;L^{p}(\mathbb{R}^n).
\end{equation}
To prove (\ref{s01}), we only need to show that for any $\epsilon >0$, there exists $R>0$ such that
$$
\int_{\mathbb{R}^n \setminus B_{R}}|u_{n}(x)|^pdx < \epsilon.$$ 
Set
\begin{eqnarray*}
&& B_{R} = \{x\in \mathbb{R}^n:\;\;|x| < R\},\\
&& A(R,M) = \{x\in \mathbb{R}^n\setminus B_{R}:\;\;V(x) \geq M\},\\
&& B(R,M) = \{x\in \mathbb{R}^n \setminus B_{R}:\;\; V(x) <M\},
\end{eqnarray*}
then 
$$
\int_{A(R, M)} |u_{n}(x)|^p dx \leq \int_{\mathbb{R}^n} \frac{V(x)}{M} |u_{n}(x)|^pdx \leq \frac{\|u_{n}\|_{X^s}^{p}}{M}.
$$
Now choose $\sigma \in (1, \frac{p_{s}^{*}}{p})$ such that
$$
\frac{1}{\sigma} + \frac{1}{\sigma'} = 1,
$$
then we have
\begin{eqnarray*}
\int_{B(R,M)} |u_{n}(x)|^pdx &\leq& \left(\int_{B(R,M)} |u_{n}(x)|^{p\sigma}dx \right)^{1/\sigma} (\mu (B(R,M)))^{1/\sigma'} \\
&\leq & C\|u_{n}\|_{X^s}^{p} (\mu (B(R,M)))^{1/\sigma'}.
\end{eqnarray*}
Since $\|u_{n}\|_{X^s}$ is bounded and condition $(V_{1})$ holds, we may choose $R, M$ large enough such that $\frac{\|u_{n}\|_{X^s}^{p}}{M}$ and $\mu (B(R,M))$ are small enough. Hence, $\forall \epsilon >0$, we have
$$
\int_{\mathbb{R}^n \setminus B_{R}} |u_{n}(x)|^p dx = \int_{A(R,M)}|u_{n}(x)|^pdx + \int_{B(R,M)} |u_{n}(x)|^pdx <\epsilon,
$$ 
from which (\ref{s01}) follows. 

To prove the lemma for general exponent $q$, we use an interpolation argument. Let $u_{n} \rightharpoonup 0$ in $X^s$, we have just proved that $u_{n} \to 0$ in $L^{p}(\mathbb{R}^n)$. That is
\begin{equation}\label{s02}
\int_{\mathbb{R}^n} |u_{n}(x)|^pdx \to 0.
\end{equation}
Moreover, because the embedding $X^s \rightharpoonup L^{p_{s}^{*}}(\mathbb{R}^n)$ is continuous and $\{u_{n}\}$ is bounded in $X^s$, we also have
\begin{equation}\label{s03}
\sup_{n} \int_{\mathbb{R}^n} |u_{n}(x)|^{p_{s}^{*}}dx < \infty.
\end{equation} 
Since $q\in (p,p_{s}^{*})$, there is a number $\lambda \in (0,1)$ such that
$$
\frac{1}{q} = \frac{\lambda}{p} + \frac{1-\lambda}{p_{s}^{*}}.
$$
Then by H\"older inequality
\begin{eqnarray*}
\int_{\mathbb{R}^n} |u_{n}(x)|^{q}dx & = & \int_{\mathbb{R}^n} |u_{n}(x)|^{\lambda q}|u_{n}(x)|^{(1-\lambda)q}dx\\
&\leq & \|u_{n}\|_{p}^{\lambda q} \|u_{n}\|_{p_{s}^{*}}^{(1-\lambda)q} \to 0.
\end{eqnarray*}
This implies $u_{n} \to 0$ in $L^{q}(\mathbb{R}^n)$. $\Box$

The dual space of $(X^s, \|.\|_{X^s})$ is denoted by $((X^s)^{*}, \|.\|_{*})$. We rephrase variationally the fractional p-Laplacian as the nonlinear operator $A: X^s \to (X^s)^{*}$  defined for all $u,v\in X^s$ by
$$
\langle A(u), v \rangle = \int_{\mathbb{R}^{n}}\int_{\mathbb{R}^{n}} \frac{|u(x) - u(y)|^{p-2}(u(x) - u(y))(v(x) - v(y))}{|x-y|^{n+sp}}dxdy + \int_{\mathbb{R}^n} V(x)|u|^{p-2}uvdx.
$$
It can be seen that, if $u$ is smooth enough, this definition coincides with that of (\ref{01}). A weak solution of problem (\ref{00}) is a function $u\in X^s$ such that
\begin{equation}\label{s03}
\langle A(u), v\rangle = \int_{\mathbb{R}^{n}} f(x,u)vdx
\end{equation}
for all $v\in X^s$. Clearly, for all $u\in X^s$
$$
\langle A(u), u\rangle = \|u\|_{X^s}^{p},\;\;\|A(u)\|_{*} \leq \|u\|_{X^s}^{p-1},
$$
and $A$ satisfies the following properties.

\begin{Lem}\label{X2lem}
For any $u,v \in X^s$, it holds that
$$
\langle A(u) - A(v), u-v \rangle \geq (\|u\|_{X^s}^{p-1} - \|v\|_{X^s}^{p-1})(\|u\|_{X^s} - \|v\|_{X^s})
$$
\end{Lem}

\noindent 
{\bf Proof.} By direct computation, we have
\begin{eqnarray*}
&&\langle A(u) - A(v), u-v \rangle =  \langle A(u), u-v \rangle - \langle A(v), u-v \rangle\\
& = &\!\!\! \int_{\mathbb{R}^{n}}\!\!\int_{\mathbb{R}^{n}} \frac{|u(x) - u(y)|^{p-2}(u(x) - u(y))((u-v)(x) - (u-v)(y))}{|x-y|^{n+sp}}dxdy +\!\! \int_{\mathbb{R}^n} V(x)|u|^{p-2}u(u-v)dx \\
&-&\!\!\! \int_{\mathbb{R}^{n}}\!\!\int_{\mathbb{R}^{n}} \frac{|v(x) - v(y)|^{p-2}(v(x) - v(y))((u-v)(x) - (u-v)(y))}{|x-y|^{n+sp}}dxdy -\!\! \int_{\mathbb{R}^n} V(x)|v|^{p-2}v(u-v)dx\\
& = & \!\!\! \int_{\mathbb{R}^{n}}\!\!\int_{\mathbb{R}^{n}} \frac{|u(x) - u(y)|^{p}}{|x-y|^{n+sp}}dxdy - \!\!\! \int_{\mathbb{R}^{n}}\!\!\int_{\mathbb{R}^{n}} \frac{|u(x) - u(y)|^{p-2}(u(x) - u(y))(v(x) - v(y))}{|x-y|^{n+sp}}dxdy\\
&&-\!\!\! \int_{\mathbb{R}^{n}}\!\!\int_{\mathbb{R}^{n}} \frac{|v(x) - v(y)|^{p-2}(v(x) - v(y))(u(x) - v(y))}{|x-y|^{n+sp}}dxdy +\!\! \int_{\mathbb{R}^{n}}\!\!\int_{\mathbb{R}^{n}} \frac{|v(x) - v(y)|^{p}}{|x-y|^{n+sp}}dxdy\\
&& +\int_{\mathbb{R}^{n}} V(x)|u|^pdx - \int_{\mathbb{R}^n}V(x)|u|^{p-2}uvdx - \int_{\mathbb{R}^n}V(x)|v|^{p-2}vudx + \int_{\mathbb{R}^n}V(x)|v|^pdx\\
& = & \|u\|_{X^s}^{p} + \|v\|_{X^s}^{p} \\
&& - \!\!\! \int_{\mathbb{R}^{n}}\!\!\int_{\mathbb{R}^{n}} \frac{|u(x) - u(y)|^{p-2}(u(x) - u(y))(v(x) - v(y))}{|x-y|^{n+sp}}dxdy - \int_{\mathbb{R}^n}\!\!V(x)|u|^{p-2}uvdx \\
&&-\!\!\! \int_{\mathbb{R}^{n}}\!\!\int_{\mathbb{R}^{n}} \frac{|v(x) - v(y)|^{p-2}(v(x) - v(y))(u(x) - v(y))}{|x-y|^{n+sp}}dxdy -\!\! \int_{\mathbb{R}^n}V(x)|v|^{p-2}vudx. 
\end{eqnarray*}
By H\"older inequality, it holds that
\begin{eqnarray*}
&&\int_{\mathbb{R}^{n}}\!\!\int_{\mathbb{R}^{n}} \frac{|u(x) - u(y)|^{p-1}|v(x) - v(y)|}{|x-y|^{n+sp}}dxdy + \int_{\mathbb{R}^n}\!\!V(x)|u|^{p-1}|v|dx\\
&\leq& \left(\int_{\mathbb{R}^n}\int_{\mathbb{R}^n} \frac{|u(x) - u(y)|^p}{|x-y|^{n+sp}}dxdy \right)^{\frac{p-1}{p}} \left(\int_{\mathbb{R}^n}\int_{\mathbb{R}^n} \frac{|v(x) - v(y)|^p}{|x-y|^{n+sp}}dxdy \right)^{\frac{1}{p}}\\
&& + \left(\int_{\mathbb{R}^n} V(x)|u|^p\right)^{\frac{p-1}{p}} \left(\int_{\mathbb{R}^n} V(x)|v|^p \right)^{\frac{1}{p}}.
\end{eqnarray*} 
Using the following inequality
$$
(a+b)^{\beta}(c+d)^{1-\beta} \geq a^{\beta} c^{1-\beta} + b^{\beta}d^{1-\beta}
$$
which holds for any $\beta \in (0,1)$ and $a>0$, $b>0$, $c>0$, $d>0$, set $\beta = \frac{p-1}{p}$ and  
\begin{eqnarray*}
a & = & \int_{\mathbb{R}^n}\int_{\mathbb{R}^n} \frac{|u(x) - u(y)|^p}{|x-y|^{n+sp}}dxdy\\
b & = & \int_{\mathbb{R}^n} V(x)|u|^pdx\\
c &= & \int_{\mathbb{R}^n}\int_{\mathbb{R}^n} \frac{|v(x) - v(y)|^p}{|x-y|^{n+sp}}dxdy\\
d & = & \int_{\mathbb{R}^n} V(x)|v|^pdx
\end{eqnarray*}
we can deduce that
$$
\int_{\mathbb{R}^{n}}\!\!\int_{\mathbb{R}^{n}} \frac{|u(x) - u(y)|^{p-2}(u(x) - u(y))(v(x) - v(y))}{|x-y|^{n+sp}}dxdy + \int_{\mathbb{R}^n}\!\!V(x)|u|^{p-2}uvdx \leq \|u\|_{X^s}^{p-1} \|v\|_{X^s} 
$$
Similarly, we can obtain 
$$
\int_{\mathbb{R}^{n}}\!\!\int_{\mathbb{R}^{n}} \frac{|v(x) - v(y)|^{p-2}(v(x) - v(y))(u(x) - u(y))}{|x-y|^{n+sp}}dxdy + \int_{\mathbb{R}^n}\!\!V(x)|v|^{p-2}uvdx \leq \|v\|_{X^s}^{p-1} \|u\|_{X^s} 
$$
Therefore we have
\begin{eqnarray*}
\langle A(u) - A(v), u-v \rangle & \geq & \|u\|_{X^s}^{p} + \|v\|_{X^s}^{p} - \|u\|_{X^s}^{p-1}\|v\|_{X^s} - \|v\|_{X^s}^{p-1}\|u\|_{X^s}\\
& = & (\|u\|_{X^s}^{p-1} - \|v\|_{X^s}^{p-1})(\|u\|_{X^s} - \|v\|_{X^s}).
\end{eqnarray*}
$\Box$

\begin{Lem}\label{X3lem}
If $u_n \rightharpoonup u$ and $\langle A(u_{n}) , u_{n} - u \rangle \to 0$, then $u_{n} \to u$ in $X^s$
\end{Lem}

\noindent
{\bf Proof.} Since $X^s$ is a reflexive Banach space, it is isometrically isomorphic to a locally uniformly convex space. So as it was proved in \cite{GD}, weak convergence and norm convergence imply strong convergence. Therefore we only need to show that $\|u_{n}\|_{X^s} \to \|u\|_{X^s}$.

We note that 
$$
\lim_{n \to \infty} \langle A(u_{n}) - A(u), u_{n} - u \rangle = \lim_{n\to \infty} (\langle A(u_{n}), u_{n} - u \rangle -  \langle A(u), u_{n} - u\rangle) =0
$$
By Lemma \ref{X2lem} we have
$$
\langle A(u_{n}) - A(u), u_n - u \rangle \geq (\|u_{n}\|_{X^s}^{p-1} - \|u\|_{X^s}^{p-1})(\|u_n\|_{X^s} - \|u\|_{X^s}) \geq 0.
$$
Hence $\|u_{n}\|_{X^s} \to \|u\|_{X^s}$ as $n\to \infty$ and the assertion follows. $\Box$





\section{Proof of Theorem \ref{Itm}}

In this section we consider the functional $I:X^s \to \mathbb{R}$ defined by
\begin{equation}\label{T01}
I(u) = \frac{1}{p}\|u\|_{X^s}^{p} - \int_{\mathbb{R}^n}F(x,u)dx.
\end{equation}
$I$ is well defined and of class $C^1$. The derivative of $I$ is given by
\begin{equation}\label{T02}
\langle I'(u), v \rangle = \langle A(u), v\rangle - \int_{\mathbb{R}^n} f(x,u)vdx,
\end{equation}
for $v\in X^s$. Therefore, the critical points  of $I$ are weak solutions of (\ref{00}).

\begin{Lem}\label{T1lem}
Suppose that $(V_{1})$, ($f_{1}$), ($f_2$) and ($f_3$) hold, then $I$ satisfies the Cerami condition (C). 
\end{Lem}

\noindent 
{\bf Proof.} Let $\{u_{k}\}$ be a sequence in $X^s$ satisfying 
$$
I(u_k) \to c,\;\;(1 + \|u_k\|_{X^s})I'(u_k) \to 0.
$$
We claim that $\{u_k\}$ is bounded in $X^s$. Otherwise, if $\|u_k\|_{X^s} \to \infty$, we consider $w_k = \frac{u_k}{\|u_k\|_{X^s}}$. Then, up to subsequence, we have
\begin{eqnarray*}
&&w_k  \rightharpoonup w\;\;\mbox{in}\;\;X^s\\
&&w_k \to w\;\;\mbox{in}\;\;L^{q}(\mathbb{R}^n)\;\;\mbox{for}\;\;p\leq q \leq p_{s}^{*}\\
&&w_k(x) \to w(x)\;\;\mbox{a.e.}\;\;x\in \mathbb{R}^n,
\end{eqnarray*}
as $k\to \infty$. We first consider the case that $w\neq 0$ in $X^s$. Since $I'(u_{k})u_k \to 0$, we have
\begin{equation}\label{T03}
\|u_k\|_{X^s}^{p} - \int_{\mathbb{R}^n} f(x,u_k)u_{k}dx \to 0.
\end{equation}
By dividing the left-hand side of (\ref{T03}) with $\|u_k\|_{X^s}^{p}$ we get
\begin{equation}\label{T04}
\left| \int_{\mathbb{R}^n} \frac{f(x,u_k)u_k}{\|u_k\|_{X^s}^{p}} \right| \leq 1.
\end{equation}
On the other hand, by Fatou's Lemma and condition $(f_1)$ we have
$$
\int_{\mathbb{R}^n} \frac{f(x,u_k)u_k}{\|u_k\|_{X^s}^{p}}dx = \int_{\{w_{k} \neq 0\}} |w_{k}|^p \frac{f(x,u_k)u_k}{|u_{k}|^p}dx \to \infty,
$$
this contradicts to (\ref{T04}).

If $w=0$ in $X^s$, consider 
\begin{eqnarray*}
\gamma_{k}:[0,1] &\to& \mathbb{R}\\
t & \to & I(tu_k).
\end{eqnarray*}
By the continuity of $\gamma_{k}$, we choose a sequence $\{t_k\} \in [0,1]$ such that
$$
I(t_{k}u_k) = \max_{t\in [0,1]}I(tu_k).
$$ 
For any $\lambda >0$, let $v_{k} = (2p\lambda)^{1/p}w_k = \frac{(2p\lambda)^{1/p}u_k}{\|u_k\|_{X^s}}$, then
\begin{equation}\label{T05}
v_k \to 0 \;\;\mbox{in}\;\;L^{q}(\mathbb{R}^n)\;\;\mbox{for}\;\;q\in [p,p_{s}^{*}).
\end{equation}
We claim that
\begin{equation}\label{T06}
\lim_{n\to \infty} \int_{\mathbb{R}^n}F(x,v_k)dx = 0.
\end{equation}
In fact, by (\ref{04}) and (\ref{T05}) for $k$ large enough we have
\begin{eqnarray*}
\left| \int_{\mathbb{R}^n}F(x,v_k)dx \right| &\leq& \frac{\epsilon}{p} \|v_{k}\|_{L^p}^{p} + \frac{C_{\epsilon}}{q}\|v_k\|_{L^q}^{q} \to 0.
\end{eqnarray*}

Since $\|u_k\|_{X^s} \to \infty$, for $k$ large enough we have
$$
\frac{(2p\lambda)^{1/p}}{\|u_k\|_{X^s}} \in (0,1).
$$ 
Hence, for $k$ large enough, (\ref{T06}) gives 
\begin{eqnarray*}
I(t_{k}u_k) & \geq & I(v_{k}) \\
& = & \frac{1}{p}\|v_k\|_{X^s}^{p} - \int_{\mathbb{R}^n}F(x,v_k)dx\\
& = & 2\lambda - \int_{\mathbb{R}^n}F(x,v_n)dx \geq \lambda.
\end{eqnarray*}
That is 
\begin{equation}\label{T07}
I(t_ku_k) \to +\infty.
\end{equation} 
Since $I(0) = 0$, $I(u_k) \to c$, we have $t_k \in (0,1)$. By the definition of $t_k$, 
\begin{equation}\label{T08}
\langle I'(t_ku_k), t_ku_k \rangle = 0.
\end{equation}
From (\ref{T07}), (\ref{T08}), we have
$$
I(t_ku_k) - \frac{1}{p}\langle I'(t_ku_k), t_ku_k\rangle = \int_{\mathbb{R}^n} \left(\frac{1}{p}f(x,t_ku_k)t_ku_k - F(x,t_ku_k)  \right)dx \to \infty.
$$
By ($f_{3}$), there exists $\theta \geq 1$ such that
\begin{equation}\label{T09}
\int_{\mathbb{R}^n} \left(\frac{1}{p}f(x,u_k)u_k - F(x,u_k)  \right)dx \geq \frac{1}{\theta} \int_{\mathbb{R}^n} \left(\frac{1}{p}f(x,t_ku_k)t_ku_k - F(x,t_ku_k)  \right)dx \to \infty.
\end{equation}
On the other hand,
\begin{equation}\label{T10}
\int_{\mathbb{R}^n} \left(\frac{1}{p}f(x,u_k)u_k - F(x,u_k)  \right)dx = I(u_{k}) - \frac{1}{p} \langle I'(u_k), u_k \rangle \to c.
\end{equation} 
(\ref{T09}) contradicts (\ref{T10}). Hence $\{u_k\}$ is bounded in $X^s$, therefore up to a subsequence we may assume that $u_k  \rightharpoonup u$ in $X^s$ and 
\begin{equation}\label{T11}
u_k \to u \;\;\mbox{in}\;\;L^{q}(\mathbb{R}^n).
\end{equation} 
By the boundedness of $\{u_k\}$ in $L^{p}(\mathbb{R}^n)$, we have
\begin{equation}\label{T12}
\Lambda_{1} = \sup_{k} \int_{\mathbb{R}^n} |u_k|^{p}dx < \infty.
\end{equation}
By H\"older inequality we also have
\begin{eqnarray*}
\int_{\mathbb{R}^n} |u_k|^{p-1}u dx & \leq & 2 \left( \int_{\mathbb{R}^n} |u_k|^{(p-1)p'}dx \right)^{1/p'} \|u\|_{L^p}\\
& = & 2\|u\|_{L^p} \left( \int_{\mathbb{R}^n} |u_k|^pdx\right)^{1/p'}.
\end{eqnarray*}
Applying (\ref{T12}), we deduce
$$
\Lambda_{2} = \sup_{k} \int_{\mathbb{R}^n} |u_k|^{p-1}udx < \infty.
$$ 
Similarly,
$$
\Lambda_3 = \sup_{k} \int_{\mathbb{R}^n} |u|^{p-1}u_{k}dx <\infty.
$$
Now, by ($f_1$) and ($f_2$), for any $\epsilon>0$ there exists $C_{\epsilon} >0$ such that
$$
|f(x,t)| \leq \epsilon |t|^{p-1} + C_{\epsilon}|t|^{q-1},\;\;\mbox{for all}\;(x,t) \in \mathbb{R}^n \times \mathbb{R}.
$$
Then using H\"older inequality we have
\begin{eqnarray}\label{T13}
&&\int_{\mathbb{R}^n} (f(x,u_k) - f(x,u))(u_{k} - u)dx \nonumber \\
&&\leq \int_{\mathbb{R}^n} [\epsilon (|u_k|^{p-1} + |u|^{p-1}) + C_{\epsilon}(|u_k|^{q-1} + |u|^{q-1})]|u_k - u|dx\nonumber\\
&&\leq \epsilon \int_{\mathbb{R}^n} (|u_k|^{p} + |u|^p + |u_{k}|^{p-1}|u| + |u|^{p-1}|u_{k}|)dx\nonumber\\
&&+C_{\epsilon} \left( \int_{\mathbb{R}^n} |u_k|^{q-1}|u_k - u|dx  + \int_{\mathbb{R}^n} |u|^{q-1}|u_k - u|dx\right)\nonumber\\
&&\leq \epsilon \left( \Lambda_{1} + \Lambda_2 + \Lambda_3 + \int_{\mathbb{R}^n} |u|^{p}dx \right)\nonumber \\
&&+ 2C_{\epsilon}\left( \sup_{k}\||u_k|^{q-1}\|_{L^{q'}} + \||u|^{q-1}\|_{L^{q'}} \right)\|u_k - u\|_{L^q}
\end{eqnarray}
Since $\{u_k\}$ is bounded in $L^q(\mathbb{R}^n)$, it follows that $\{|u_k|^{q-1}\}$ is bounded in $L^{q'}(\mathbb{R}^n)$. That is 
$$
\sup_{k} \||u_k|^{q-1}\|_{L^{q'}} < \infty.
$$
Therefore we can deduce from (\ref{T11}) and (\ref{T13}) that
$$
\int_{\mathbb{R}^n} (f(x,u_k) - f(x,u))(u_k - u)dx \to 0.
$$
Note that $I'(u_k) \to 0$ and we have
$$
\langle A(u_k) - A(u), u_k - u \rangle = \langle I'(u_k) - I'(u), u_k - u \rangle + \int_{\mathbb{R}^n} (f(x,u_k) - f(x,u))(u_k - u)dx \to 0.
$$
Therefore, by Lemma \ref{X3lem} we obtain $u_n \to u$ in $X^s$. The proof is complete. $\Box$

\noindent
{\bf Proof of Theorem \ref{Itm}}

We check that $I$ has the mountain pass geometry. By (\ref{04}), we have
$$
|F(x,t)| \leq \frac{\epsilon}{p} |t|^p + \frac{C_\epsilon}{q}|t|^q,
$$
so
\begin{eqnarray*}
I(u) &\geq& \frac{1}{p} \|u\|_{X^s}^{p} - \frac{\epsilon}{p} \|u\|_{L^p}^{p} - \frac{C_\epsilon}{q}\|u\|_{L^q}^{q}\\
&\geq& \left( \frac{1}{p} - \frac{\epsilon C}{p}  \right)\|u\|_{X^s}^{p} - \frac{C_\epsilon C}{q}\|u\|_{X^s}^{q}
\end{eqnarray*}
Let $\epsilon>0$ small enough such that $\frac{1}{p} - \frac{\epsilon C}{p} >0$ and $\|u\|_{X^s} = \rho$. Since $q>p$, taking $\rho$ small enough such that
$$
\frac{1}{p} - \frac{\epsilon C}{p} - \frac{C_\epsilon C}{q}\rho^{q-p} >0.
$$
Therefore 
$$
I(u) \geq \rho^{p} \left( \frac{1}{p} - \frac{\epsilon C}{p} - \frac{C_\epsilon C}{q}\rho^{q-p} \right) = \beta >0.
$$

Now we note that by ($f_1$)
$$
\lim_{|t| \to \infty} \frac{f(x,t)t}{|t|^p} = +\infty\;\;\mbox{implies that} \;\;\lim_{|t|\to \infty} \frac{F(x,t)}{|t|^p} = +\infty.
$$
So, for any $\epsilon >0$, there exists $M>0$ such that
$$
F(x,t) > \frac{|t|^p}{\epsilon},\;\;\mbox{for all}\;|t| >M.
$$
Let $c(\epsilon) = \frac{M^p}{\epsilon}$, then
$$
F(x,t) > \frac{|t|^p}{\epsilon} - \frac{M^p}{\epsilon}.
$$
Next, for $\varphi \in C_{0}^{\infty}(\mathbb{R}^n)$ we have
$$
\int_{\mathbb{R}^n} \frac{F(x,t\varphi)}{|t|^p} \geq \frac{1}{\epsilon}\int_{\mathbb{R}^n} |\varphi|^pdx - \frac{M^p}{\epsilon |t|^p} \int_{supp (\varphi)}dx. 
$$
This implies
\begin{equation}\label{T14}
\lim_{|t| \to \infty} \int_{\mathbb{R}^n} \frac{F(x,t\varphi)}{|t|^p}dx \geq \frac{1}{\epsilon}\int_{\mathbb{R}^n} |\varphi|^pdx,
\end{equation}
for all $\epsilon >0$. Since $\epsilon$ is arbitrary, by (\ref{T14}) we get
$$
\lim_{|t| \to \infty} \int_{\mathbb{R}^n} \frac{F(x,t\varphi)}{|t|^p}dx = +\infty.
$$ 
Consequently,
$$
\frac{I(t\varphi)}{|t|^p} = \frac{1}{p}\|\varphi\|_{X^s}^{p} - \int_{\mathbb{R}^n} \frac{F(x,t\varphi)}{|t|^p}dx \to -\infty,\;\;\mbox{as}\;|t| \to \infty.
$$ 
Hence, let $t_{0}$ be big enough and $e = t_{0}\varphi$, then we have $I(e) < 0$.

Therefore, since by Lemma \ref{T1lem}, $I$ satisfies the Cerami condition and has mountain pass geometry, using the Mountain pass Lemma the proof of Theorem is complete. $\Box$

\section{Symmetry Results}


In order to prove Theorem \ref{SRtm1}, our main tools will be symmetry rearrangement. For the reader convenience we remember some result about symmetric decreasing rearrangement. For more details see \cite{ELML}.  

Let $A\subset \mathbb{R}^{n}$ be a Lebesgue measurable set and denote the measure of $A$ by $|A|$. Define the symmetrization $A^{*}$ of $A$ to be the closed ball centered at the origin such with the same measure as $A$. Thus, if we define $\omega (n)$ to be the volume of the unit ball in $\mathbb{R}^{n}$, then for $A \subset \mathbb{R}^{n}$
$$
A^{*} := B(0, (|A|/\omega (n))^{1/n}).
$$

Let $u:\mathbb{R}^{n} \to \mathbb{R}$ a Borel measurable function, then $u$ is said to vanish at infinity if
$$
|\{x: |u(x)|>t\}|< \infty\;\;\mbox{for all}\;\;t>0.
$$

Now if $u:\mathbb{R}^{n} \to \mathbb{R}$ is a Borel measurable function vanishing at infinity we define
\begin{equation}\label{S01}
u^{*}(x) = \int_{0}^{\infty} \chi_{\{|u|\geq t\}}^{*}(x)dt,
\end{equation}
where $\chi_{A}^{*}:= \chi_{A^{*}}$.

The rearrangement $u^{*}$ has a number of properties, see \cite{ELML}:
\begin{itemize}
\item[(i)] $u^{*}$ is nonnegative.
\item[(ii)] $u^{*}$ is radially symmetric and nonincreasing, i.e:
\begin{eqnarray*}
|x| \leq |y| \;\;\mbox{implies}\;\;u^{*}(y) \leq u^{*}(x)
\end{eqnarray*}
\item[(iii)] $u^{*}$ is a lower semicontinuos function.
\item[(iv)] The level sets of $u^{*}$ are simply the rearrangement of the level set of $u$, i.e
 $$
 \{x:u^{*}(x) > t\} = \{x:|u(x)|>t\}^{*}.
 $$
 an important consequence of this is the equimeasurability  of the function $u$ and $u^{*}$, i.e
\begin{eqnarray*}
|\{u^{*} > t\}| = |\{|u| > t\}|\;\;\mbox{for all}\;\;t>0.
\end{eqnarray*}
\item[(v)] For any positive monotone function $\phi$, we have
$$
\int_{\mathbb{R}^{n}} \phi (|u(x)|)dx = \int_{\mathbb{R}^{n}} \phi (u^{*}(x))dx.
$$
In particular, $u^{*}\in L^{p}(\mathbb{R}^{n})$ if and only if $u\in L^{p}(\mathbb{R}^{n})$ and
$$
\|u\|_{L^{p}} = \|u^{*}\|_{L^{p}}
$$
\item[(vi)] Let $V(|x|) \geq 0$ be a spherically symmetric increasing function on $\mathbb{R}^{n}$. If $u$ is a nonnegative function on $\mathbb{R}^{n}$, vanishing at infinity the
    $$
    \int_{\mathbb{R}^{n}} V(|x|)|u^{*}(x)|^{2}dx \leq \int_{\mathbb{R}^{n}}V(|x|)|u(x)|^{2}dx
    $$
\item[(vii)] \textbf{Riesz' rearrangement inequality.} Let $u, v, w$ be nonnegative measurable functions on $\mathbb{R}^{n}$ that vanish at infinity. Then
    $$
    \int_{\mathbb{R}^{n}}\int_{\mathbb{R}^{n}}u(x)v(x-y)w(y)dydx \leq \int_{\mathbb{R}^{n}}\int_{\mathbb{R}^{n}}u^{*}(x)v^{*}(x-y)w^{*}(y)dydx
    $$
\end{itemize}

\begin{Thm}\label{RPStm2}
{\cite{FAEL} \bf p-Polya-Szeg\"o inequality.} Let $0< s < 1$ and  $1\leq p < \infty$. Define $u^* = |u|^*$, then 
\begin{enumerate}
\item IF $u\in W^{s, p}(\mathbb{R}^n)$, then $u^*$ is also in $W^{s , p}(\mathbb{R}^n)$ and
\begin{equation}\label{S03}
\int_{\mathbb{R}^n}\int_{\mathbb{R}^n} \frac{|u^{*}(x) - u^{*}(z)|^p}{|x-z|^{n+sp}}dzdx \leq \int_{\mathbb{R}^n}\int_{\mathbb{R}^n} \frac{|u(x) - u(z)|^p}{|x-z|^{n+sp}}dzdx.
\end{equation}
\item If $u_{k}$ is a sequence of functions in $W^{s,p}(\mathbb{R}^n)$ and if $u_k \to u$ strongly in $W^{s,p}(\mathbb{R}^n)$ norm, then $u_k^* \to u^*$ strongly in $W^{s,p}(\mathbb{R}^n)$ norm.
\end{enumerate}
\end{Thm}

\noindent 
{\bf Proof of Theorem \ref{SRtm1}}

By the proof of Theorem \ref{Itm}, under ($f_{1}$)-($f_{3}$) and ($V$)  we find that $I$ satisfies the mountain pass geometry conditions and the mountain pass level for $I$ is given by
$$
c= \inf_{\gamma \in \Gamma}\sup_{t\in [0,1]}I(\gamma (t)),
$$
where $\Gamma$ is defined as usual.
By definition of $
c$, for any $k\in \mathbb{N}$, there is $\gamma_{k}\in \Gamma$ such that
\begin{equation}\label{SReq11}
\sup_{t\in[0,1]}I(\gamma_{k}(t)) \leq c + \frac{1}{k^{2}}.
\end{equation}
By the triangle inequality 
$$
||u_k(x+z)| - |u_k(x)|| \leq |u_k(x+z) - u_k(x)|,
$$
thus the norm $\|.\|_{X^s}$ of $|u_k|$ is not bigger that the one of $u_k$. Therefore, without loss of generality, we may suppose that $u_k$ is nonnegative.
Now, let $\gamma_{k}^{*}(t) = [\gamma_{k}(t)]^{*}$. By the continuity of rearrangements  in $X^{s}$ we have that $\gamma_{k}^{*} \in \Gamma$. Moreover, by Theorem \ref{RPStm2} and taking into account that $V$ satisfies ($V_r$), we have
$$
I(\gamma_{k}^{*}(t)) \leq I(\gamma_{k}(t)), \;\;\forall t\in [0,1].
$$
So
\begin{equation}\label{SReq12}
\sup_{t\in [0,1]}I(\gamma_{k}^{*}(t)) \leq c + \frac{1}{k^{2}}.
\end{equation}
By the Ekeland's variational principle \cite{MW}, there is a sequence $w_{k} \in X^s$ and $\xi_{k}\in [0,1]$ such that
\begin{equation}\label{SReq13}
\| w_{k} - \gamma_{k}^{*}(\xi_{k}) \|_{X^s} \leq \frac{1}{k},
\end{equation}
\begin{equation}\label{SReq14}
I(w_{k}) \in (c - \frac{1}{k^{2}}, c + \frac{1}{k^{2}}),
\end{equation}
\begin{equation}\label{SReq15}
\| I'(w_{k}) \|_{(X^s)'} \leq \frac{1}{k}.
\end{equation}
Following the ideas of the proof of Theorem \ref{Itm}, we can show that
$ w_{k} \to  w,$ $I(w) =  c$, $I'(w)w  =  0$ and finally that
\begin{equation}\label{SReq16}
\lim_{k\to \infty}\|w - \gamma_{k}^{*}(\xi_{k})\|_{X^s} = 0,
\end{equation}
concluding the proof. $\Box$


\medskip


\end{document}